# MINIMUM MULTIPLICITIES OF SUBGRAPHS AND HAMILTONIAN CYCLES

**JOHN SHEEHAN**

Department of Mathematical Sciences
University of Aberdeen
King's College
Aberdeen AB24 3UE

# Abstract

Let $G$ be a finite simple graph with automorphism group $A(G)$. Then a spanning subgraph $U$ of $G$ is a *fixing subgraph* of $G$ if $G$ contains exactly $|A(G)|/|A(G) \cap A(U)|$ subgraphs isomorphic to $U$: the graph $G$ must always contain at least this number. If in addition $A(U) \subseteq A(G)$ then $U$ is a *strong fixing subgraph*. Fixing subgraphs are important in many areas of graph theory. We consider them in the context of Hamiltonian graphs.

## §0  Summary of contents

Let $G$ be a finite simple graph with automorphism group $A(G)$. Suppose that $U$ is a spanning subgraph of $G$. Then $U$ is a *fixing subgraph of G* if $G$ contains exactly $|A(G)|/|A(G) \cap A(U)|$ subgraphs isomorphic to $U$: if in addition $A(U) \subseteq A(G)$ then $U$ is said to be a *strong fixing subgraph of G*. In a sense fixing subgraphs have minimum multiplicity, up to isomorphism, in $G$.

The idea of a *strong* fixing subgraph (originally simply called a fixing subgraph) was introduced in [12-13]. In this paper we consider a less restrictive definition of a fixing subgraph. The underlying idea is of course classical (see for example Fraïssé [8] and Droste [7]). Where our approach and the classical approach diverge is that the interest here is in the local properties of finite graphs rather than the global properties of infinite graphs.

We believe that fixing subgraphs are of interest in many areas of graph theory: in this paper their significance is considered in the context of Hamiltonian graphs. We consider fixing subgraphs generally in section 1, in the general context of hamiltonian cycles in section 2 and specifically generalized Petersen graphs in section 3.

The reason for the name "fixing" hopefully should become more obvious in section 1: if $U$ is a strong fixing subgraph of $G$ it "fixes" the number of extensions of $U$ to $G$.

## §1  Introduction

Let $G$ be a finite simple graph with vertex set $V(G)$, edge set $E(G)$ and automorphism group $A(G)$. Let $[G]$ denote the isomorphism class of graphs containing $G$.

Let span($G$) denote the set of spanning subgraphs of $G$: if $U \in$ span($G$) write $U \subseteq G$.

Suppose that $U \subseteq G$. Set $[U;G] = [U] \cap$ span($G$) and set

$$s(U;G) = \|[U;G]\|,$$

that is, $s(U;G)$ is the number of subgraphs of $G$ which are isomorphic to $U$.

Finally set

$$x(U;G) = |\{X : U \subseteq X \cong G\}|,$$

that is, $x(U;G)$ is the number of "extensions" of $U$ to a graph $X$ isomorphic to $G$.

**Theorem 1** [4]

$$|A(U)|s(U;G) = |A(G)|x(U;G).$$

**Proof**  Let $\mathcal{S}_n$ be the symmetric group acting on $V(G)$ (where $|V(G)| = n$) with induced action on pairs of vertices defined by $\pi(u,v) = (\pi(u), \pi(v))$. Then $\pi(G)$ is the graph with vertex set $V(G)$ and edge-set $\{\pi(e) : e \in E(G)\}$. We have

$$|A(U)|s(U;G) = |\{\pi \in \mathcal{S}_n : \pi(U) \subseteq G\}|$$
$$= |\{\sigma \in \mathcal{S}_n : U \subseteq \sigma(G)\}|$$
$$= |A(G)|x(U;G). \quad \bullet$$

**Notation**  Suppose that $U, U' \in \text{span}(G)$. Then $U$ and $U'$ are said to be *G-similar* if there exists $\sigma \in A(G)$ such that $U' = \sigma(U)$. Set

$$s_0(U;G) = |\{U' : U' \subseteq G, U' \text{ is } G\text{-similar to } U\}|.$$

**Theorem 2**

$$s(U;G) \geq s_0(U;G) = \frac{|A(G)|}{|A(U) \cap A(G)|} \quad (1)$$

$$\geq \frac{|A(G)|}{|A(U)|}. \quad (2)$$

**Proof**  Suppose that $U_i$ is $G$-similar to $U$ ($i = 1,2$). Then there exist $\sigma_i \in A(G)$ such that $\sigma_i(U_i) = U$. Hence $\sigma_1(U) = \sigma_2(U)$ and $\sigma_1^{-1}\sigma_2 \in A(G) \cap A(U)$. Hence $s_0(U;G) = |A(G)|/|A(U) \cap A(G)|$. This is the only part of the statement which is non-trivial.  $\bullet$

**MAIN DEFINITION**   Suppose $U \in \text{span}(G)$. Then $U$ is said to be a *fixing subgraph* of $G$ if equality holds throughout (1): $U$ is said to be a *strong fixing subgraph* of $G$ if equality holds in both (1) and (2).  •

The classical definition for general graphs (see [8]) is given in Theorem 3 which is simply a restatement of Theorem 2.

**Theorem 3**   Let $U \in \text{span}(G)$. Then

(i)   $U$ is a fixing subgraph of $G$ if and only if for each $U_0 \in [U;G]$ there exists an isomorphism $\sigma : U \to U_0$ such that $\sigma$ extends to an automorphism of $G$.

(ii)   $U$ is a strong fixing subgraph of $G$ if and only if for each $U_0 \in [U,G]$ and each isomorphism $\sigma : U_0 \to U$, $\sigma$ extends to an automorphism of $G$.

**Proof**   Suppose that $U \in F(G)$ and $U_0 \in [U;G]$. Then, from Theorem 2, $U_0$ is $G$-similar to $U$.  •

**Corollary 4**   Suppose that $U$ is a fixing subgraph of $G$ and $U_0 \in [U;G]$. Then $U_0$ is a fixing subgraph of $G$ and if $U$ is a strong fixing subgraph then so is $U_0$.
•

Set
$$F(G) = \{[U;G] : U \text{ is a fixing subgraph of } G\}$$
and
$$F^*(G) = \{[U;G] : U \text{ is a strong fixing subgraph of } G\}.$$
To avoid repetition below the superscript '$*$' is always used to indicate the analogous definition for "strong" fixing subgraphs. Unless any ambiguity arises we simply ignore the distinction between a graph and its isomorphism type, for example, we write $U \in F(G)$ rather than $[U;G] \subseteq F(G)$.

The next theorem gives some simple consequences of the definitions.

**Theorem 5**   Suppose that $U \in \text{span}(G)$. Then

(i)   $U \in F(G)$ if and only if $x(U;G) = \dfrac{|A(U)|}{|A(U) \cap A(G)|}$.

(ii)   $U \in F^*(G)$ if and only if $x(U;G) = 1$.

(iii)   Suppose that $U \in F^*(G)$ and $U \subseteq K \subseteq G$. Then $A(K) \subseteq A(G)$ and $K \in F^*(G)$.

**Proof**   (i) and (ii)   This follows immediately from the definition of a (strong) fixing subgraph and Theorem 1.

(iii)   Suppose that $U \in F^*(G)$ and $U \subseteq K \subseteq G$. Then, from (ii), $x(U;G) = 1$ and hence $x(K;G) = 1$. Therefore $K \in F^*(G)$ and by definition $A(K) \subseteq A(G)$.   □

**Comment**   Theorem 5(iii) underlies the reason why the definition of strong fixing subgraphs is so restrictive. Their structure is more accessible of course: in a sense the minimal elements of $F^*(G)$ determine $F^*(G)$ whereas this is clearly not the case with $F(G)$.   □

In the next section we illustrate these ideas in the context of hamiltonian cycles.

## §2   Hamiltonian cycles

To test whether a spanning subgraph $U$ is a (strong) fixing subgraph of $G$, there are three useful methods: (i) the direct approach of Theorem 3; (ii) estimating $x(U;G)$ and (iii) estimating $s(U;G)$. We demonstrate these methods in the next example.

**Example**   Let $C$ be the hamiltonian cycle

$$C = (u_0, u_1, u_2, ..., u_{13})$$

of the Heawood graph $H$. We do not distinguish here, and elsewhere, between the cycle $C$ and the spanning subgraph whose edges are the edges of the cycle. Integers are modulo 14.

(i) The mapping $\sigma : u_i \to u_{i+1}$ $(i = 0,\ldots,13)$ is an automorphism of $C$ but $\sigma \notin A(G)$. Hence, from Theorem 3, $C \notin F^*(G)$.

(ii) There are exactly two 'extensions' of $C$ to $H$, that is, $x(C;H) = 2$. This is easy to verify since $H$ is bipartite and the *girth* $\gamma(H)$ of $H$ is 6. Thus starting with the cycle $C$ consider which vertex is adjacent to $u_0$ other than $u_1$ and $u_{13}$. Suppose that $u_0 \sim u_7$ (read as $u_0$ is adjacent to $u_7$). Then $u_{13} \sim u_4$ and $u_{12}$ is adjacent only to $u_{11}$ and $u_{13}$ which is impossible. Therefore $u_0 \sim u_5$ or $u_0 \sim u_9$. Without loss of generality choose $u_0 \sim u_5$ then $u_1 \sim u_{10}$ and so on. Since $|A(C)| = 28$ and $|A(C) \cap A(H)| = 14$, $x(C;H) = |A(C)|/|A(C) \cap A(H)|$ and $C \in F(G)$.

(iii) We prove again that $C$ is a fixing subgraph by examining $s(C;H)$. One way of checking that $H$ contains exactly 24 hamiltonian cycles is as follows. Let $C_0$ be any hamiltonian cycle of $H$ and define $p = p(C_0,C)$ to be the length of a longest common path of $C_0$ and $C$, for example, if

$$C_0 = (u_0, u_1, u_{10}, u_{11}, u_6, u_7, u_2, u_3, u_{12}, u_{13}, u_8, u_9, u_4, u_5)$$

then $p(C_0,C) = 1$. By inspection, there exist 2, 7, 7, 7, 1 hamiltonian cycles corresponding respectively to $p$ equal to 1, 2, 3, 4, 14. Now, since $H$ is cubic, 4-transitive and not 5-transitive $|A(H)| = 14 \times 3 \times 2^3$ and therefore $s(C;H) = |A(H)|/|A(C) \cap A(H)|$. □

Let $\mathcal{F}(HAM)$ be the set of Hamiltonian graphs $G$ such that for each hamiltonian cycle $C$ of $G$, $C \in F(G)$.

**Theorem 6** Suppose that $G$ has $n$ ($\geq 3$) vertices. Then $G \in \mathcal{F}^*(HAM)$ if and only if $G$ is isomorphic to either: (i) $C_n$; (ii) $K_n$ or (iii) $K_{n/2,n/2}$ (when $n$ is even). □

**Comment** (i) The proof is given in the appendix: in fact a stronger result is proved since only the local property that the automorphisms of the hamiltonian cycle(s) extend to automorphisms of $G$, is used.

(ii) Theorem 6 is somewhat disappointing. This prompted an investigation of $\mathcal{F}(\text{HAM})$ where the structure is much richer. We are a long way from any worthwhile characterization. Trivially if $G$ has a unique hamiltonian cycle [2, 15-17] then $G \in \mathcal{F}(\text{HAM})$. The author [15] incidentally conjectured in 1975 that every 4-regular Hamiltonian graph has at least two hamiltonian cycles. This in turn would imply that, apart from the cycle, every regular Hamiltonian graph has at least two hamiltonian cycles.

From the example above the Heawood graph belongs to $\mathcal{F}(\text{HAM})$: this graph is often [18] called Tutte's 6-cage (possibly iterated star products of the Heawood graph also belong to $\mathcal{F}(\text{HAM})$). It is easy to check (see also Theorem 6) that the 3-cage $K_4$ and the 4-cage $K_{3,3}$ belong to $\mathcal{F}^*(\text{HAM})$. Not so easily we have proved, without the use of a computer, that Tutte's 8-cage belongs to $\mathcal{F}(\text{HAM})$. This graph $G$ has 144 hamiltonian cycles: it is 5-transitive with $|A(G)| = 30 \times 3 \times 2^4$ and $|A(G) \cap A(C)| = 10$ where $C$ is any hamiltonian cycle.

## §3  Generalized Petersen graphs and hamiltonian cycles

The hamiltonian structure of Generalized Petersen graphs has been much studied and the classification [1] of the Hamiltonian Generalized Petersen graphs proved difficult.

The Generalized Petersen graph $G = G(n,k)$, $1 \leq k < n/2$, is defined by
$$V(G) = \{a_i : i = 1,2,...,n\} \cup \{b_i : i = 1,2,...,n\}$$
$$E(G) = \{a_i a_{i+1} : i = 1,2,...,n\} \cup \{a_i b_i : i = 1,2,...,n\}$$
$$\cup \{b_i b_{i+k} : i = 1,2,...,n\}$$

(where integers are taken modulo $n$).

The cycle $(a_1, a_2, \ldots, a_n)$ is called the *rim* of $G$.

Bondy [3] gives a detailed proof that $G(n,2)$ is non-hamiltonian if and only if $n \equiv 5 \pmod 6$. We can prove that:

**Theorem 7** (i) $G(n,2) \in \mathcal{F}(\text{HAM})$ if and only if $n \equiv 1,3 \pmod 6$ or $n = 10$.

(ii) $G(n,1) \in \mathcal{F}(\text{HAM})$ if and only if $n$ is odd or $n = 4$.

**Proof** (i) Set $G = G(n,2)$. Define permutations $\alpha$ and $\beta$ on $V(G)$ by

$$\alpha(a_i) = a_{i+1}, \quad \alpha(b_i) = b_{i+1}$$
$$\beta(a_i) = a_{-i}, \quad \beta(b_i) = b_{-i}$$

for $i = 1, 2, \ldots, n$. Let $D_n = \langle \alpha, \beta \rangle$, that is, the dihedral group of order $2n$.

**Claim** [9]   $A(G) = D_n$  $(n \neq 5, 10)$.

If $n \equiv 5 \pmod 6$ then $G$ is non-hamiltonian and so $G \notin \mathcal{F}(\text{HAM})$.

So now suppose $n \not\equiv 5 \pmod 6$ and $n \neq 10$. Then $G$ is hamiltonian. Consider a hamiltonian cycle of $G$. The intersection of this cycle with the rim of $G$ will be a sequence of paths of lengths $n_1, n_2, \ldots, n_r$ taken in order round the rim. By specifying this sequence and using the fact that a hamiltonian cycle is a connected 2-factor it is routine to reconstruct possible hamiltonian cycles.

Using the claim it is also routine to test whether two given sequences correspond to $G$-similar hamiltonian cycles.

In each case below we identify the sequence $(n_1, n_2, \ldots, n_r)$ with its associated hamiltonian cycle(s).

**Case 1**   $n \equiv 5 \pmod 6$.

In this case $G$ is not hamiltonian and so $G \notin \mathcal{F}(\text{HAM})$.

**Case 2**   $n \equiv 0, 2, 4 \pmod 6$  $(n \neq 10)$.

Choose any $k$, $1 \leq k < n$, so that $k \equiv n - 3 \pmod 4$ and $k \not\equiv 0 \pmod 3$ except possibly $k = 3$. Then $G$ contains the hamiltonian cycle $(n_1, n_2, \ldots, n_r)$ where $n_1 = k$ and $n_i = 1$ $(i > 1)$. From the claim distinct $k$'s give hamiltonian cycles which are not $G$-similar. Hence $G \notin \mathcal{F}(\text{HAM})$.

**Case 3**   $n \equiv 1, 3 \pmod{6}$.

From [3, p.61] since $n$ is odd, $n_i \leq 2$ for each $i$.

When $n \equiv 1 \pmod{6}$ the only sequence with $n_1 \geq n_2 \geq \ldots \geq n_r$ which is hamiltonian has $n_i = 2$ $(i = 1, \ldots, r-2)$ and $n_{r-1} = n_{r-2} = 1$.

When $n \equiv 3 \pmod{6}$ the only sequence which is hamiltonian has $n_i = 2$ for each $i$.

In each case, from the claim, the cycles are $G$-similar and $G \in \mathcal{F}(\text{HAM})$.

**Case 4**   $(n = 10)$

$G = G(10,2)$ is the dodecahedron. By inspection $G$ contains exactly 30 hamiltonian cycles $C$, $|A(G)| = 120$ and $|A(C) \cap A(G)| = 4$. Hence
$$s(C; G) = |A(G)| / |A(C) \cap A(G)|$$
and $G \in \mathcal{F}(\text{HAM})$.  □

(ii) Set $G = G(n,1)$. $G$ contains the hamiltonian cycle
$$((n-1)), \tag{3}$$
that is, $n_1 = n - 1$ and $r = 1$.

When $n$ is odd this is the only type of hamiltonian cycle contained in $G$ and clearly all such cycles are $G$-similar. Hence $G \in \mathcal{F}(\text{HAM})$.

When $n$ is even $G$ contains, in addition, the hamiltonian cycle
$$(1,1,\ldots,1), \tag{4}$$
that is, $n_i = 1$ $(i = 1, 2, \ldots, n/2)$.

The cycles of types (3) and (4) are not $G$-similar $(n > 4)$: the non-edges of the cycle of type (4) are all 4-chords which is not the case with the cycle of type (3). Hence $G \notin \mathcal{F}(\text{HAM})$. When $n = 4$ these types of cycles are $G$-similar and hence the cube, $G(4,1) \in \mathcal{F}(\text{HAM})$.  □

In [3] Bondy states that $G(n,3)$ is hamiltonian. In proving Theorem 8 we verify this statement. The main details are given in the appendix.

**Theorem 8**  $G(n,3) \notin \mathcal{F}(\text{HAM})$  $(n \neq 8, 10)$.  □

**Comment**  Suppose that $(n,k) \in \{(4,1),(8,3),(10,2),(10,3),(12,5),(24,5)\}$. Then the graphs $G(n,k)$ are the "very symmetric" exceptional graphs noted in [9, Theorem 2]. In a sense the higher the symmetry of $G$ the more likely it is that $G \in \mathcal{F}(\text{HAM})$. On the other hand, the homogeneous [5-6, 10-11, 14] graph, $L(K_{3,3})$, which is highly symmetric, does not belong to $\mathcal{F}(\text{HAM})$.  □

## §4.  Appendix

**Proof**  (Theorem 6)

Let $G \in \mathcal{F}^*(\text{HAM})$. Suppose that $G$ contains the hamiltonian cycle

$$C = (u_0, u_1, \ldots, u_{n-1}) \qquad (n \geq 2)$$

where integers are now, usually, considered modulo $n$. We do not distinguish between the cycle $C$ and the spanning subgraph whose edges are the edges of the cycle. Since $C \subseteq G$, $A(C) \subseteq A(G)$. Therefore there exists $D \subseteq \mathbb{Z}_n$ such that

$$u_i \sim u_j \text{ if and only if } i - j = D$$

(notice that since $A(C)$ is the dihedral group and $A(C) \subseteq A(G)$, $i - j \in D$ if and only if $j - i \in D$). Hence, setting $|D| = k$, $G$ is a circulant graph of degree $k$.

If $k = 2$ then $G \cong C_n$ and the theorem is true. So now assume that $k \geq 3$. Choose $m$ so that $u_0 \sim u_m$, $2 \leq m \leq \left\lfloor \frac{n}{2} \right\rfloor$ and $m$ is as large as possible (because $G$ is a circulant this slight abuse of terminology - there being no ordering in $\mathbb{Z}_n$ - is harmless). Since $m \in D$, $u_1 \sim u_{m+1}$.

Therefore $G$ contains the hamiltonian cycle

$$C_0 = (u_0, u_m, u_{m-1}, \ldots, u_1, u_{m+1}, u_{m+2}, \ldots, u_{n-1})$$

and again $A(C_0) \subseteq A(G)$. Therefore

$$u_{m-i} \sim u_{m+i+1} \quad (i = 0,1,2,...,m\text{-}1),$$

that is, $2i+1 \in D$ $(i = 0,1,...,m\text{-}1)$. By the maximality of $m$ and since $i-j \in D$ if and only if $j-i \in D$, a contradiction is avoided only if $G \cong K_n$ when $n$ is odd or $K_{n/2,n/2} \subseteq G$ when $n$ is even. So now assume that $n$ $(\geq 4)$ is even and $G \not\cong K_{n/2,n/2}$.

$$\text{Then } 2i-1 \in D \quad (i = 1,2,...,n/2): \text{ in particular} \tag{A.1}$$
$$u_0 \sim u_{2i-1} \quad (i = 1,2,...,n/2).$$

Now choose $m$ so that $u_0 \sim u_m$, $m$ is even, $n/2 \leq m \leq n-2$ and $m$ is as large as possible. Since $m \in D$, $u_1 \sim u_{m+1}$. Set $m = 2s$. Then $G$ contains the hamiltonian cycle

$$C_1 = (u_0, u_{2s}, u_{2s-1}, ..., u_2, u_1, u_{2s+1}, u_{2s+2}, ..., u_{n-1}).$$

Then, replacing $C$ by $C_1$ in (A.1), $u_0 \sim u_{2k}$ $(k = 1,2,...,s)$. It follows from the maximality of $m$ and since $i-j \in D$ if and only if $j-i \in D$, that $m = n-2$ and $G \cong K_n$.

Sufficiency is easy to prove. □

**Proof** (Theorem 8)

Set $G = G(n,3)$ $(n \neq 8,10)$.

**Claim [9]** $A(G) = D_n$.

As in the proof of Theorem 7 we simply list, in each case, two sequences $(n_1, n_2, ..., n_r)$ (identified with hamiltonian cycles in $G$) which, using the claim, are not $G$-similar.

**Case 1** ($n$ even)

$n_i = 1$ $(i = 1,2,...,n/2)$.

**Subcase 1.1** ($n \equiv 0 \pmod 4$)

(a) $n \equiv 0 \pmod 3$; set $t = (n-6)/3$:

$n_i = 2$ $(i = 1,2,...,t)$;  $n_j = 1$ $(j = t+1, t+2, t+3)$.

(b) $n \not\equiv 0 \pmod{3}$; set $t = n/4$:

$n_i = 3$ $(i = 1,2,...,t)$.

**Subcase 1.2** ($n \equiv 2 \pmod 4$)

Set $t = (n-6)/4$.

(a) $n \equiv 0 \pmod 3$;

$n_i = 3$ $(i = 1,2,...,t)$;  $n_j = 1$ $(j = t+1, t+2, t+3)$.

(b) $n \not\equiv 0 \pmod 3$;

$n_i = 3$ $(i = 1,2,...,t)$;  $n_j = 1$ $(j = t+1, t+2, t+3)$.

**Case 2** ($n$ odd)

(a) $n \equiv 2 \pmod 3$;

(i) $n_1 = n - 13$; $n_2 = n_3 = 1$; $n_4 = 4$; $n_5 = 2$

(ii) Set $t = (n-8)/3$:

$n_i = 2$ $(i = 1,2,...,t)$;  $n_j = 1$ $(j = t+1,...,t+4)$.

(b) $n \not\equiv 2 \pmod 3$;

(i) $n_1 = n - 5$; $n_2 = n_3 = 1$

(ii) $n_1 = n - 11$; $n_2 = n_3 = 2$; $n_4 = n_5 = 1$.  □